\documentclass[lettersize,journal]{IEEEtran}
\usepackage{amsmath,amsfonts}
\usepackage{algorithmic}
\usepackage{algorithm}
\usepackage{array}
\usepackage[caption=false,font=normalsize,labelfont=sf,textfont=sf]{subfig}
\usepackage{textcomp}
\usepackage{stfloats}
\usepackage{url,color}
\usepackage{verbatim}
\usepackage{graphicx}
\usepackage{cite,bm}

\usepackage{enumerate}
\usepackage{latexsym,amsmath,amssymb,amsfonts,epsfig,graphicx,psfrag}

\newtheorem{lemma}{Lemma}
\newtheorem{theorem}{Theorem}

\newtheorem{assumption}{Assumption}

\begin{document}

\title{Distributed Optimization with Projection-free Dynamics}

\author{
	Guanpu Chen, Peng Yi, Yiguang Hong~\IEEEmembership{Fellow,~IEEE,} Jie Chen~\IEEEmembership{Fellow,~IEEE,}
\thanks{This work was supported by the National Natural Science Foundation of China (No. 61733018),
        	and by Shanghai Municipal Science and Technology Major Project (No. 2021SHZDZX0100).}
\thanks{G. Chen is with the Key Laboratory of Systems and Control, Academy of Mathematics and Systems Science, Beijing, China.
	{\tt\small chengp@amss.ac.cn}}%
\thanks{P. Yi is with Department of Control Science and Engineering \& Shanghai Research Institute for Intelligent Autonomous Systems, Tongji University, Shanghai, China.
	{\tt\small yipeng@tongji.edu.cn}}%
\thanks{Y. Hong is with Department of Control Science and Engineering \& Shanghai Research Institute for Intelligent Autonomous Systems, Tongji University, Shanghai, China, and is also with the Key Laboratory of Systems and Control, Academy of Mathematics and Systems Science, Beijing, China.
	{\tt\small yghong@iss.ac.cn}}
\thanks{J. Chen is with Department of Control Science and Engineering \& Shanghai Research Institute for Intelligent Autonomous Systems, Tongji University, Shanghai, China.
	{\tt\small chenjie206@tongji.edu.cn}}
}



\maketitle

\begin{abstract}
We consider continuous-time dynamics for distributed optimization with set constraints in the paper. To handle the computational complexity of projection-based dynamics due to solving a general quadratic optimization subproblem with projection, we propose a distributed projection-free dynamics by employing the Frank-Wolfe method, also known as the conditional gradient algorithm. The process searches a feasible descent direction by solving an alternative linear optimization instead of a quadratic one. To make the approach applicable over weight-balanced digraphs, we design a dynamics for the consensus of local decision variables and another dynamics of auxiliary variables to track the global gradient. Then we prove the convergence of the dynamical systems to the optimal solution, and provide detailed numerical comparisons with both projection-based dynamics and other distributed projection-free algorithms. Also, we derive the distributed discrete-time scheme following the instructive ideas of the proposed dynamics and provide its accordingly convergence rate.
\end{abstract}


\section{Introduction}
Distributed optimization and its applications have attracted a large amount of research attention in the past decade. Under multi-agent frameworks, the global objective function consists of agents' local objective functions, and each agent shares limited amounts of information with neighbors through the networks to achieve an optimal solution. Both discrete-time algorithms \cite{nedic2010constrained,yuan2012distributed,li2020distributed,yuan2015regularized} and continuous-time algorithms \cite{wang2015distributed,yang2016multi,zeng2016distributed,yi2018distributed,shi2018finite,lu2018distributed} are extensively developed for solving distributed problems.

Among continuous-time algorithms, projection-based dynamics have been widely adopted to solve distributed optimization with constraints, on the basis of the well-developed theory in nonlinear optimization \cite{khalil2002nonlinear,ruszczynski2011nonlinear}.
Various projection-based dynamics have been designed with techniques in dynamical systems and control theory.
For example,
 \cite{yang2016multi} proposed a proportional-integral protocol to solve distributed constrained optimization, while \cite{zeng2016distributed} proposed  distributed dynamics where the projection maps are with respect to tangent cones.
However, projection-based design implies that agents will encounter a quadratic subproblem when a variable needs to find its nearest point to a set. When the constraints are expressed as complex structures, the computational cost of quadratic subproblems discourages agents from employing projection-based approaches, particularly for high-dimensional constraints.

Meanwhile, though the primal-dual method is another common tool for distributed constrained optimization \cite{yuan2015regularized,yi2016initialization,zhu2018projected,yang2019survey,liang2020distributed}, the attendant difficulties and challenges actually occur in the high dimension of the algorithm, since the dimension of dual dynamics is as high as the dimension of the constraints. This brings extra burden in local storage, algorithm iteration and distributed communication. Furthermore, practical constraints are usually equipped with sparsity, that is, we may spend huge calculation resources to iterate the additional dual dynamics. Nowadays, considering the optimization with large-scale data or constraints, such as sparse PCA,  structural SVM, and matrix completion, researchers are eagerly seeking efficient and specific tools for optimization with high-dimensional constraints, rather than the common approaches like projection-based or primal-dual algorithms.

Fortunately, the well-known Frank-Wolfe (FW) method \cite{frank1956algorithm}, also known as the conditional gradient, provides us with a helpful idea. Briefly speaking, the FW method uses a linearized function to approximate the objective function and derives a feasible descent direction by solving optimization with a linear objective. Thanks to the linear programming toolbox, the  feasible descent direction can be efficiently computed when the constraints are high-dimensional polyhedrons, which can usually be used as an approximation for general convex sets \cite{chen123distributed}. Then, this process avoids general projection operations and solving corresponding quadratic subproblems in algorithm iterations.
	In fact, the FW method is suitable for various sparse optimization with atomic constraints \cite{jaggi2013revisiting}, among which the linear constraint is the most comprehensible one. Though proposed early in 1956 \cite{frank1956algorithm}, the FW method has regained lots of attention in recent years due to its computational advantage for optimization  with  large-scale sparse constraints, which covers wide applications such as image process in machine learning and data analysis in bioinformation \cite{cui2014generalized,chapel2020partial,zhang2017bio,jiang2022distributed}.
	Also, there have been massive investigations on the FW method afterward, such as rate analysis over strongly convex sets in \cite{garber2015faster}, online learning over networks in \cite{zhang2017projection},
	and quantized FW for lower communication in \cite{zhang2020quantized}. 
%

Motivated by the above, we aim to design a projection-free dynamics based on the FW method for solving distributed constrained optimization in this paper. 
Agents have their own local objective functions and need to achieve the optimal solution via communicating with neighbors over networks. The main contributions are as follows.
 \subsection*{Contributions}

 Inspired by the historical literature \cite{jacimovic1999continuous}, we design a novel FW-based distributed dynamics for agents to solve the constrained optimization. Compared to the dynamics in \cite{yang2016multi,zeng2016distributed}, a feasible descent direction is derived by solving  optimization with a linear objective. Hence, the dynamics avoids solving complicated quadratic subproblems due to projection operations on set constraints, which actually leads to a projection-free dynamics. 
Our distributed projection-free dynamics accordingly induces a novel discrete-time format, where averaging consensus is employed both  to ensure the consensus of local decision variables and to  help auxiliary variables  track the global gradient. {This differs from the mechanism in the decentralized discrete-time FW algorithm of \cite{wai2017decentralized}. Hence, we develop a novel convergence analysis with the Lyapunov theory and comparison theorems, and provide the analysis on the convergence rate.
 Moreover, compared with the projected dynamics given in \cite{yang2016multi,zeng2016distributed} and the  discrete-time FW algorithm in \cite{wai2017decentralized}, the distributed projection-free dynamics is applied over the communication networks described by weight-balanced digraphs.
	
\subsection*{Paper Organization} 
The organization  of this paper is given as follows.
Section \ref{Formulation and Algorithm} formulates the distributed constrained optimization with basic assumptions, and presents a  projection-free algorithm with the corresponding discretization.
Then Section \ref{main} reports the main results, including the consensus of decision variables, global gradient tracking and convergence of the continuous-time algorithm, as well as the convergence rate of the induced discrete scheme. Finally, Section \ref{example} shows numerical examples with a comparison to the existing algorithms, while Section \ref{conclusion} provides some concluding remarks and discussions on future work.

\subsection*{Notations} 
	Here we give some necessary notations in this paper.
	Denote $\mathbb R^n$ ( $\mathbb R^{m\times n}$) as the set of $n$-dimension ($m$-by-$n$) real column vectors (real matrices), $ 1_{n}$ (or $ 0_{n}$) as the $n$-dimensional column vectors with all elements of one (or zero), and
	$col\{x_1,\dots,x_n\}=(x_1^{\rm T},\dots,x_n^{\rm T})^{\rm T}.$
	Let $A \otimes B$ as the Kronecker product of matrices $A$ and $B$. In addition, take $\|\cdot\|$  as the Euclidean norm of vectors, and	$\|\|_F$ as the Frobenius norm of real matrices defined by 
	$\|Q\|_F=\sqrt{tr(Q^TQ)}.$

\section{Distributed Projection-free Dynamics}\label{Formulation and Algorithm}

In this section, we formulate the constrained distributed optimization, and then propose the distributed projection-free dynamics, as well as a discrete scheme derived from the continuous-time algorithm.
\subsection{Problem Formulation}
Consider $N$ agents indexed by $\mathcal V=\{1,2\dots,N\}$. For agent $i\in\mathcal V$, there is  a local cost function $f_i:\mathbb R^n\rightarrow \mathbb R$ on the feasible set $\Omega\subseteq \mathbb R^n$. The global cost function is
$$F(x)=\frac{1}{N}\sum_{i=1}^N f_i(x).$$
All agents aim to solve the  constrained optimization:
\begin{align}\label{ori_problem}
\min_{x\in\mathbb R^n}F(x)\quad \text{s.t.}, \; x\in\Omega.
\end{align}

In a multi-agent network, the $i$th agent controls a local decision variable $x_i\in\Omega$ to search the optimal $x^* \in \arg\min F(x)\quad \text{s.t.}, \; x\in\Omega.$ Also, the information of local cost functions are regarded as private knowledge.
The agents communicate with their neighbors through a network described by a digraph $\mathcal G(\mathcal V, \mathcal E)$, where $\mathcal V$ is the set of nodes (regarded as agents here) and $\mathcal E\subseteq\mathcal V \times \mathcal V$ is the set of edges. $\mathcal A=[a_{ij}]\in\mathbb R^{N\times N}$ is the adjacency matrix subject to $a_{ij}>0$ if and only if  $(i,j)\in\mathcal E$, which means that agent $j$ can send information to agent $i$, and $a_{ij}=0$, otherwise.  The  Laplacian matrix  $\mathcal L=\mathcal D- \mathcal A$, where $\mathcal D\in\mathbb R^{N\times N}$ is diagonal with $\mathcal D_{i,i}=\sum_{j=1}^N a_{ij}$, for any $i\in\mathcal V$.
A digraph $\mathcal G$ is strongly connected if there exists at least one directed path between any pair of vertices, and $\mathcal G$ is weight-balanced if $\sum_{j=1}^N a_{ij}=\sum_{j=1}^N a_{ji}$ for $i\in\mathcal V$.

In addition,
the gradient of a differentiable function $f$ is $\kappa$-Lipschitz on convex set $C\subseteq \mathbb R^n$ with a constant $\kappa>0$, if
$$\|\nabla f(x)-\nabla f(y)\|\leq \kappa \|x-y\|,\quad\forall x,y \in C.$$
Also, the above is equivalent to the following:
$$f(x)-f(y)\leq(x-y)^{\rm T}\nabla f(y)+\frac{\kappa}{2}\|x-y\|^2,\quad\forall x,y \in C.$$

Then we consider solving distributed optimization \eqref{ori_problem} under
the following assumptions.
\begin{assumption}\label{assum}
	\
	\begin{itemize}
		\item The feasible set $\Omega$ is convex, compact and nonempty.
		\item For $i\in\mathcal V$, $f_i$ is convex and differentiable, and $\nabla f_i$ is $\kappa$-Lipschitz on $\Omega$.
		\item The digraph $\mathcal G$ is strongly connected and weight-balanced.
	\end{itemize}	
\end{assumption}
The differentiable cost functions enable us to use the gradient information as in \cite{garber2015faster,wai2017decentralized,zhang2017projection}. Additionally, the strongly connected and weight-balanced digraph, as a generalization of connected undirected graphs, was studied in other continuous-time distributed algorithms  \cite{gharesifard2013distributed,liang2020distributed}.
\subsection{Distributed Projection-free Dynamics}
To solve the distributed optimization \eqref{ori_problem}, 
we intend to explore a projection-free approach following the FW method. 
In fact, continuous-time algorithms provide evolved dynamics, which help reveal
how each variable achieves the optimal solution with given directions. Also, the analysis techniques in modern calculus and nonlinear systems may lead to comprehensive results with mild assumptions in continuous-time schemes. Hence, we firstly propose a novel projection-free dynamics with the FW method in the following \textit{Algorithm \ref{alg:1}}, which differs from the projection-based continuous-time algorithms in  \cite{yang2016multi,zeng2016distributed}.
\begin{algorithm}[H]
	\caption{Distributed Projection-free Dynamics for $i\in\mathcal V$}
	\label{alg:1}
\vspace{5pt}
	\textbf{Initialization:}
	
\vspace{5pt}
	$x_i(0)\in\Omega$, $y_i(0)=\mathbf{0}_n$, $z_i(0)\in\mathbb R^n$, $v_i(0)\in\Omega$.
	
\vspace{5pt}
	\textbf{Flows renewal:}
	\begin{flalign*}
	\dot x_i(t) =&\sum_{j=1}^{N}a_{ij}(x_j(t)-x_i(t))+\beta(t)(v_i(t)-x_i(t)),&\\
		v_i(t)  \in&~\arg\min_{v\in\Omega}~ z_i(t)^Tv,&\\
	\dot y_i(t) =&\sum_{j=1}^{N}a_{ij}(z_j(t)-z_i(t)),&\\
	z_i(t) =&~ y_i(t)+\nabla f_i(x_i(t)),& 
	\end{flalign*}
	where $\beta(t)$ is a positive time-varying parameter satisfying \begin{flalign*}
&\lim_{t\rightarrow \infty}\beta(t)=0,\quad \lim_{t\rightarrow \infty}\int_{0}^{t}\beta(\tau)d\tau=\infty.&
	\end{flalign*}
\end{algorithm}
\textit{Algorithm \ref{alg:1}} is distributed since the dynamics of the $i$th agent is only influenced by local values of $x_i$, $y_i$, $z_i$, $v_i$ and $\nabla f_i(x_i)$. Specifically, the $i$th agent uses local decision variable $x_i$ to estimate the optimal solution $x^*\in\Omega$
and  local optimal solution $v_i$ as a conditional gradient. Since each agent is merely capable to calculate its own gradient $\nabla f_i(x)$, rather than the global gradient {$\frac{1}{N}\sum_{i=1}^{N}\nabla f_i(x)$}, thus,  $z_i$ serves as the variable that simultaneously operates two processes --- one is to compute agent $i$'s local gradient, the other is to achieve consensus with neighbors' local gradients, in order to estimate the global gradient. In fact, the gradient tracking method in \cite{nedic2017achieving,pu2020distributed} motivates our algorithm design. Although the time-varying $\beta(t)$ seems to be a global parameter, it is easy to determine its value for all agents, by merely selecting some general decreasing functions like $\beta(t)=1/t$.  That is analogous to other FW-based works dealing with parameters \cite{wai2017decentralized,zhang2017projection}.
\subsection{Discretization}
So far, we have devoted to proposing distributed projection-free dynamics to avoid solving complex subproblems owing to projections.
Following these instructive ideas in designing the above distributed continuous-time Algorithm 1, we consider deriving its implementable discretization. 
We notice the decentralized discrete-time FW algorithm in \cite{wai2017decentralized}.
To make a comparison, we present the discretization of \textit{Algorithm \ref{alg:1}} by adopting the notations in \cite{wai2017decentralized} in the following. To remain consistent with  \cite{wai2017decentralized}, set the network $\mathcal G$ undirected and connected, and adjacency matrix $\mathcal A$ as symmetric and doubly stochastic.
Let $0<\delta<1$ be a fixed step-size in discretization, and denote $\eta^k=\delta\beta^k$.
Take the weighted average of agent $i$'s neighbors $\mathcal N_i$ in the network $\mathcal G$ by $$Avg_{j\in\mathcal N_{i}}\{x_j^k\}=(1-\delta)x_i^k+\delta \sum_{j=1}^N a_{ij}x_j^k,$$  and
similarly, 
$$Avg_{j\in\mathcal N_{i}}\{z_j^k\}=(1-\delta)z_i^k+\delta \sum_{j=1}^N a_{ij}z_j^k.$$ 
 
Then, consider the ODE in Algorithm \ref{alg:1}
\begin{align*}
\dot z_i(t) =\sum_{j=1}^{N}a_{ij}(z_j(t)-z_i(t))+\frac{d}{dt}\nabla f_i(x_i(t)).
\end{align*}
The corresponding difference equation is
\begin{align*}
z_i^{k+1} =(1-\delta)z_i^k+\delta \sum_{j=1}^N a_{ij}z_j^k+\nabla f_i(x_i^{k+1})-\nabla f_i(x_i^{k}).
\end{align*}
Therefore, the discretization of \textit{Algorithm \ref{alg:1}} gives
\begin{equation}\label{discrete-1}
\vspace{10pt}
\left\{
\begin{aligned}
x_i^{k+1}=&Avg_{j\in\mathcal N_{i}}\{x_j^k\}+\eta^k(v_i^k-x_i^k),\\
z_i^{k+1} =&Avg_{j\in\mathcal N_{i}}\{ z_j^k\} +\nabla f_i(x_i^{k+1})-\nabla f_i(x_i^k),\\
v_i^k\in&\arg\min_{v\in\Omega}~ v^Tz_i^{k}.
\end{aligned}
\right.
\end{equation}
For clarification, the discrete-time FW  in \cite{wai2017decentralized} is given as
\begin{equation}\label{discrete-2}
\vspace{10pt}
\left\{
\begin{aligned}
x_i^{k+1}=&Avg_{j\in\mathcal N_i}\{x_j^k\}+\eta^k(v_i^k-Avg_{j\in\mathcal N_i}\{x_j^k\}),\\
z_i^{k+1} =&Avg_{j\in\mathcal N_i}\{\nabla f_j(Avg_{l\in\mathcal N_j}\{ x_l^k\})\} ,\\
v_i^k\in&\arg\min_{v\in\Omega}~ v^Tz_i^{k}.
\end{aligned}
\right.
\end{equation}

The  discretization above reveals that the major difference between \eqref{discrete-1} and \eqref{discrete-2} refers to the update protocol of $z_i^{k+1}$. In \eqref{discrete-1}, agent $i$ uses both its neighbors' gradient values and its own gradient renewal to track the global gradient. In \eqref{discrete-2}, agent $i$  gathers the average value of neighbors' decision variables to compute local gradient at  first. Then agent $i$ makes again the neighbors' average gradient value to estimate the global gradient. 
In other words, the consensus format in \eqref{discrete-2} requires agents to \textit{sequentially collect} its neighbors' decision variables $x_{j\in\mathcal N_i}^k$ and auxiliary variables $z_{j\in\mathcal N_i}^k$, which results in more communication burden and data storage. Consider that constructing a reliable communication link is time-consuming and energy-consuming in practice. 
Therefore, \eqref{discrete-1} overcomes this weak point by requiring agents to \textit{simultaneously collect} the both variables.
Since the mechanism of \textit{{Algorithm \ref{alg:1}}}  differs from what in \cite{wai2017decentralized}, we need to explore new tools for analysis.
In addition, we apply \textit{{Algorithm \ref{alg:1}}} over weight-balanced digraphs as the generalization of undirected connected graphs in \cite{wai2017decentralized}, which brings more technical challenges correspondingly.
\section{Convergence Analysis}\label{main}
In this section, we give the convergence analysis of both \textit{Algorithm \ref{alg:1}} and its discretazition \eqref{discrete-1}. 
\subsection{Convergence of Algorithm \ref{alg:1}}
For convenience of statements in the following, we first provide some compact notations here. Take $$\bm\Omega\triangleq\Omega\times\cdots\times\Omega,$$
$$\bm x\triangleq col\{x_1,\dots,x_N\},$$ 
$$\bm y\triangleq col\{y_1,\dots,y_N\},$$  
$$\bm z\triangleq col\{z_1,\dots,z_N\},$$ and $\bm L=\mathcal L\otimes I_n$.
Also, define the profile of gradients
$$G(\bm x)\triangleq col\{\nabla f_1(x_1),\dots,\nabla f_N(x_N)\}.$$ 
Equivalently, \textit{Algorithm \ref{alg:1}} can be expressed in the following compact form
\begin{equation}\label{alg__exact_comp} \left\{
\begin{aligned}
\dot {\bm x}(t)=&~-\bm L\bm x(t)+\beta(t)(\bm{v}(t)-{\bm x}(t)),\\
\dot {\bm y}(t) =&~-\bm L\bm z(t),\\
{\bm z}(t) =&~{\bm y}(t)+G(\bm x(t)),
\end{aligned}
\right.
\end{equation}
where $\bm v\triangleq col\{v_1,\dots,v_N\}$ with 
$$v_i  \in \arg\min_{v\in\Omega}~ z_i^Tv.$$
Hereupon, we can derive the following results for the convergence guarantees of Algorithm \ref{alg:1}.

\begin{theorem}\label{thm_1}
	Under \textit{Assumption \ref{assum}} and with given initial condition $x_i(0)\in\Omega$, $y_i(0)=\bm 0_n$, $z_i(0)\in\mathbb R^n$ and $v_i(0)\in\Omega$,
	\begin{enumerate}[i).]
		\item all decision variables $x_i$ achieve consensus, \textit{i.e.},
		$$\lim_{t\rightarrow \infty}(x_i(t)-x_j(t))=\bm 0_n,\quad\forall i,j\in\mathcal V;$$
		\item all variables $z_i$ asymptotically track the dynamical global gradient, \textit{i.e.},
		$$\lim_{t\rightarrow \infty}\big(z_i(t)-\frac{1}{N}\sum_{j=1}^{N}\nabla f_j(x_j(t))\big)=\bm 0_n,\quad\forall i\in\mathcal V;$$
		\item 	 all decision variables $x_i$, for $i\in \mathcal{V}$, converge to a consensual optimal solution to problem \eqref{ori_problem}.
	\end{enumerate}
\end{theorem}

The following two  lemmas are necessary for the analysis of \textit{Algorithm \ref{alg:1}}, whose proofs can be found in the appendix.
\begin{lemma}\label{lemma1}
	Under \textit{Assumption \ref{assum}}, if $x_i(0)\in\Omega$ for all $i\in\mathcal V$, then $x_i(t)\in\Omega$ for all $t>0$ and for all $i\in \mathcal{V}$, .
\end{lemma}
\begin{lemma}\label{lemma2}
	Given  scalars $\varepsilon(t)\geq 0$, $s(t)\geq 0$, and  $\gamma(t)> 0$, if $\lim_{t\rightarrow \infty}\int_{0}^{t}\gamma(\tau)d\tau=\infty$, $\lim_{t\rightarrow \infty}\varepsilon(t)=0$, and
	$$\dot s(t)\leq -\gamma(t)s(t)+\gamma(t)\varepsilon(t),$$
	then $\lim_{t\rightarrow \infty}s(t)=0$.
\end{lemma}

In this way,  we give the proof of  \textit{Theorem \ref{thm_1}}.

\noindent\textbf{Proof.}

i).	It follows from \textit{Lemma \ref{lemma1}} that $x_i(t)\in{{\Omega}}$. Moreover, since $v_i(t)$ is chosen from $\Omega$, it implies that  $v_i(t)-x_i(t)$ is bounded. Because $\beta(t)\rightarrow 0$ as $t\rightarrow \infty$, we have $$\beta(t)(v_i(t)-x_i(t))\rightarrow \bm 0_n,\quad\text{as}\quad t\rightarrow \infty.$$  Thus, the dynamics for decision variables in \textit{Algorithm \ref{alg:1}} can be written as
\begin{align}\label{consensus}
\dot x_i(t) =\sum_{j=1}^{N}a_{ij}(x_j(t)-x_i(t))+u_i(t),
\end{align}
where $\lim_{t\rightarrow \infty}u_i(t)=\bm 0$.
According to the existing results in \cite[Section 4.3]{qiu2016distributed}, all decision variables in \eqref{consensus} reach consensus, \textit{i.e.}, $$\lim_{t\rightarrow \infty}(x_i(t)-x_j(t))=\bm 0_n,\quad \forall\; i,j\in\mathcal V,$$
which yields the result.

\noindent ii). Set $\hat{\bm 1}=1_N\otimes I_n$, and let us investigate
$$W(t)=\bm z(t)-\frac{1}{N} \hat{\bm 1}\hat{\bm 1}^TG(\bm x(t)).$$
Considering the orthogonal decomposition in the subspace $\ker(\bm L)$ and its complementary space $\ker(\bm L)_{\perp}$, for $t\geq 0$, define
$$W(t)=W_0(t)+W_{\perp},$$ 
and
$$\bm z(t)=\bm z_0(t)+\bm z_{\perp}(t),$$
where 
$$W_0,\bm z_0\in\ker(\bm L)=\text{span}\{ 1_N\otimes v:v\in\mathbb R^n\},$$ and $W_{\perp},\bm z_{\perp}\in\ker(\bm L)_{\perp}$. Since $ \hat{\bm 1}\hat{\bm 1}^TG(\bm x(t))\in\ker(\bm L)$, clearly, we further obtain that
$$W_0(t)=\bm z_0-\frac{1}{N} \hat{\bm 1}\hat{\bm 1}^TG(\bm x(t)),$$
and
$$ W_{\perp}(t)=\bm z_{\perp}(t),\quad \forall t\geq 0.$$
Due to the  weight-balanced digraph $\mathcal G$,
$$\sum_{i=1}^N\dot y_i(t) =\sum_{i=1}^{N}\sum_{j=1}^{N}a_{ij}(z_j(t)-z_i(t))=\bm 0_n,\quad \forall t\geq 0.$$
Together with the initial condition  $y_i(0)=\bm 0_n$ for $i\in\mathcal V$, we can verify that 
$$\sum_{i=1}^N y_i(t) =\bm 0_n.$$
Hence,
\begin{align*}\sum_{i=1}^N z_i(t)&=\sum_{i=1}^N y_i(t)+\sum_{i=1}^N\nabla f_i(x_i(t))\\
&=\sum_{i=1}^N\nabla f_i(x_i(t)),
\end{align*}
or equivalently,
$$\sum_{i=1}^N z_i(t)=\hat{\bm 1}^TG(\bm x(t)).$$
It follows from $\bm z_0\in\ker(\bm L)$ that $z_{i0}(t)=z_{j0}(t)$. Therefore, 
$$W_0(t)=\bm z_0-\frac{1}{N} \hat{\bm 1}\hat{\bm 1}^TG(\bm x(t))=\bm 0_{nN},$$
which indicates that 
$$W(t)=W_{\perp}(t).$$


%

Futhermore, let us invetigate the  function $$J(t)=\frac{1}{2}\|W(t)\|^2,$$ and consider its derivative along Algorithm \ref{alg:1} in the following.
\begin{align*}
\dot J(t)=&\Big(\bm z(t)-\frac{1}{N}\hat{\bm 1} \hat{\bm 1}^TG(\bm x(t))\Big)^T\Big(\dot{\bm z}(t)-\frac{1}{N}\hat{\bm 1} \hat{\bm 1}^T\dot G(\bm x(t))\Big)\\
=&\Big(-\bm z(t)+\frac{1}{N}\hat{\bm 1} \hat{\bm 1}^TG(\bm x(t))\Big)^T\bm L\bm z(t)\\
&+\Big(\bm z(t)-\frac{1}{N}\hat{\bm 1}\hat{\bm 1}^TG(\bm x(t))\Big)^T\Big(\bm I-\frac{1}{N}\hat{\bm 1} \hat{\bm 1}^T\Big)\dot G(\bm x(t)),
\end{align*}
where $\bm I=I_{nN}$.
By \textit{Assumption \ref{assum}}, the digraph is strongly connected and weight-balanced, which yields $$\bm L^T\hat{\bm 1}=\bm 0_{nN}.$$ Hence,
\begin{align*}
\dot J(t)=& -W(t)^T\bm LW(t)+W(t)^T\Big(\bm I-\frac{1}{N}\hat{\bm 1}\hat{\bm 1}^T\Big)\dot G(\bm x(t))\\
\leq& -W(t)^T\Big(\frac{1}{2}(\bm L+\bm L^T)\Big)W(t)\\
&+\|W(t)\|\|\bm I-\frac{1}{N}\hat{\bm 1}\hat{\bm 1}^T\|_F\|\dot G(\bm x(t))\|\\
\leq& -\lambda_2\|W(t)\|^2+\|W(t)\|\|\bm I-\frac{1}{N}\hat{\bm 1}\hat{\bm 1}^T\|_F\|\dot G(\bm x(t))\|,
\end{align*}
where $\lambda_2$ is the smallest positive eigenvalue of $\frac{1}{2}(\bm L+\bm L^T)$, and the last inequality follows from the fact $W(t)=W_{\perp}(t)$ and Rayleigh quotient theorem \cite[Page 234]{horn2012matrix}.
Moreover,
\begin{align}\label{ineq_1}
\frac{d}{dt}\|W(t)\|=&\frac{d}{dt} \sqrt{2J(t)}\nonumber\\
	=&\frac{\dot J(t)}{\|W(t)\|}\\
\leq&-\lambda_2\|W(t)\|+\|\bm I-\frac{1}{N}\hat{\bm 1}\hat{\bm 1}^T\|_F\|\dot G(\bm x(t))\|.\nonumber
\end{align}
It follows from \textit{Assumption \ref{assum}} that $\nabla f_i$ is $\kappa$-Lipschitz on $\Omega$, which leads to the $\frac{k}{2}$ boundedness of $\|\nabla G(x)\|$. Thus,
\begin{align*}
\|\dot G(\bm x(t))\|\leq& \frac{\kappa}{2}\|\dot {\bm x}(t)\|\\ =&\frac{\kappa}{2}\|\bm L\bm x(t)+\beta(t)(\bm v(t) -\bm x(t))\|.
\end{align*}
So far, since $x_i(t)$ achieves consensus and  $\beta(t)\rightarrow 0 $ as $t\rightarrow\infty$, we have
$$\|\bm L\bm x(t)+\beta(t)(\bm v(t) -\bm x(t))\|\rightarrow 0,\quad\text{as}\quad t\rightarrow \infty,$$  
which indicates that $$\lim_{t\rightarrow \infty}\|\dot G(\bm x(t))\|=0.$$	
Recall the properties that $\lambda_2>0$ and $\|\bm I-\frac{1}{N}\hat{\bm 1}\hat{\bm 1}^T\|_F>0$. Therefore,  we learn from \textit{Lemma \ref{lemma2}} that \eqref{ineq_1} yields 
$$\lim_{t\rightarrow \infty}\|W(t)\|=0,$$ which implies that all variables $z_i$ track the dynamical global gradient.

\noindent iii). Suppose that $x^*$ is an optimal solution to problem \eqref{ori_problem} and denote $\bm x^*=col\{x^*,\dots,x^*\}$. Take 
$$\overline x (t)=\frac{1}{N}\sum_{i=1}^{N}x_i(t),\quad \overline{\bm x}=col\{\overline x,\dots,\overline x\}.$$ Consider the following function
\begin{align*}
V(t)=&F(\overline x(t))-F(x^*)\\
=&\frac{1}{N}\sum_{i=1}^N f_i(\overline x (t))-\frac{1}{N}\sum_{i=1}^N f_i(x^*).
\end{align*}
Clearly, $V(t)\geq 0$. Then we investigate its derivative.
\begin{align*}
\dot V(t)=&\frac{1}{N}\Big(\sum_{i=1}^N\nabla  f_i(\overline x (t))\Big)^T\dot{\overline x} (t)\\
=&\frac{1}{N^2}\Big(\hat{\bm 1}^TG(\overline {\bm x}(t))\Big)^T\hat{\bm 1}^T \dot{\bm x}\\
=&\frac{1}{N^2}G^T(\overline {\bm x}(t))\hat{\bm 1}\hat{\bm 1}^T\Big(-\bm L\bm x(t)+\beta(t)(\bm v(t)-\bm x(t))\Big)\\
=& \frac{\beta(t)}{N^2}G^T(\overline {\bm x}(t))\hat{\bm 1}\hat{\bm 1}^T\Big(\bm v(t)-\bm x(t)\Big),
\end{align*}
where the last equality holds since $\hat{\bm 1}^T\bm L=\bm 0$.
Then,
\begin{align*}
\dot V(t)=& \frac{\beta(t)}{N^2}\Big(G^T(\overline {\bm x}(t))\hat{\bm 1}\hat{\bm 1}^T-N\bm z^T(t)\Big)^T\Big(\bm v(t)-\bm x(t)\Big)\\
&+\frac{\beta(t)}{N}\bm z^T(t)\Big(\bm v(t)-\bm x(t)\Big).
\end{align*}
Recall the derivation of $\bm  v(t)$, or equivalently, $$v_i(t)=\arg\min_{v\in\Omega} z_i^T(t)v,$$
which implies that, 
$$\bm z^T(t)\bm v(t)\leq z^T(t)\bm x', \quad \forall \bm x'\in\bm \Omega.$$ 
On this basis, take $\bm x'=\bm x^*$, and thus,
\begin{align*}
\dot V(t)\leq&~ \frac{\beta(t)}{N^2}\Big(G^T(\overline {\bm x}(t))\hat{\bm 1}\hat{\bm 1}^T-N\bm z^T(t)\Big)\Big(\bm v(t)-\bm x(t)\Big)\\
&+\frac{\beta(t)}{N}\bm z^T(t)\Big(\bm x^*-\bm x(t)\Big)\\
=&~\underbrace{\frac{\beta(t)}{N^2}G^T(\overline {\bm x}(t))\hat{\bm 1}\hat{\bm 1}^T\Big(\bm x^*-\bm x(t)\Big)}_{U_1(t)}\\
&+\underbrace{\frac{\beta(t)}{N^2}\Big(G^T(\overline {\bm x}(t))\hat{\bm 1}\hat{\bm 1}^T-N\bm z^T(t)\Big)\Big(\bm v(t)-\bm x^*\Big)}_{U_2(t)}.
\end{align*}
By the convexity of the cost functions,
\begin{align*}
U_1(t)=&\beta(t)\Big(\frac{1}{N}\sum_{i=1}^N \nabla f_i(\overline x (t))\Big)^T\Big(x^*-\overline x\Big)\\
\leq&\beta(t)\Big(-\frac{1}{N}\sum_{i=1}^N f_i(\overline x (t))+\frac{1}{N}\sum_{i=1}^N f_i(x^*)\Big)\\
=& -\beta(t)V(t).
\end{align*}
Meanwhile, it follows from  $W(t)=\bm z(t)-\frac{1}{N}\hat{\bm 1}\hat{\bm 1}^TG(\bm x(t))$ that
\begin{align*}
U_2(t)=&\frac{\beta(t)}{N^2}\Big(G^T( {\bm x}(t))\hat{\bm 1}\hat{\bm 1}^T-N\bm z^T(t)\Big)\Big(\bm v(t)-\bm x^*\Big)\\
&+\frac{\beta(t)}{N^2}\Big(G^T(\overline {\bm x}(t))-G^T( {\bm x}(t))\Big)\hat{\bm 1}\hat{\bm 1}^T\Big(\bm v(t)-\bm x^*\Big)\\
\leq&\frac{\beta(t)}{N}\Big(\|W(t)\|+\kappa\|\overline{\bm x}(t)-\bm x(t)\|\Big)\|\bm v(t)-\bm x^*\|.
\end{align*}
Since $\bm  v(t),\bm x^*\in\bm \Omega$, there exists a contant $c>0$ such that $\|\bm v(t)-\bm x^*\|\leq c$ and $$U_2(t)\leq\frac{c\beta(t)}{N}\Big(\|W(t)\|+\kappa\|\overline{\bm x}(t)-\bm x(t)\|\Big).$$

Therefore,
\begin{align*}
\dot V(t)\leq -\beta(t)V(t)+\frac{c\beta(t)}{N}\Big(\|W(t)\|+\kappa\|\overline{\bm x}(t)-\bm x(t)\|\Big).
\end{align*}
Recall that  we have already proved  $$\lim_{t\rightarrow \infty}\|\overline{\bm x}(t)-\bm x(t)\|=0,$$
as well as $$\lim_{t\rightarrow \infty}\|W(t)\|=0.$$ Moreover, the positive parameter $\beta(t)$ satisfies 
$$\lim_{t\rightarrow \infty}\int_{0}^{t}\beta(\tau)d\tau=\infty.$$ Hence, by \textit{Lemma \ref{lemma2}} again, we have  
$$\lim_{t\rightarrow \infty}V(t)=0.$$

Take $X^*\subseteq\Omega$ as the set of optimal solutions to problem \eqref{ori_problem}, and 
$$\rho (x,X^*)=\inf_{x'\in X^*}\|x-x'\|.$$ Since $\Omega$ is compact, there exists a point $x_\infty\in\Omega$ and a sequence $\{\overline x(t_k),k\in\mathbb N\}$ such that $$\lim_{k\rightarrow\infty}\|x_\infty-\overline x(t_k)\|=0,$$ and
$$\lim\sup_{t\rightarrow \infty}\rho(\overline x(t),X^*)=\lim_{k\rightarrow \infty}\rho(\overline x(t_k),X^*).$$
Since $f_i$ is differentiable  and $\rho$ is lower semicontinuous,
$$\lim_{k\rightarrow \infty}\rho(\overline x(t_k),X^*)=\rho(x_\infty,X^*).$$
Thus,
$$\lim_{k\rightarrow \infty}F(\overline x(t_k))=F(x_\infty)=F(x^*),$$
which implies $x_\infty\in X^*$, \textit{i.e.},  the decision variable $ x_i$  converges to an optimal solution to problem \eqref{ori_problem}.  \hfill$\square$	
\subsection{Convergence Rate of Discretization}
	To compare our discretization \eqref{discrete-1} with the discrete-time algorithm \eqref{discrete-2} in \cite{wai2017decentralized}, we provide theoretical guarantees on the convergence rate of discretization  \eqref{discrete-1}. Similarly, here we need some new notations to describe the discrete scheme. Take $ \boldsymbol{x}^{k}=col\{ x^{k}_{i} \}^{N}_{i=1}$, $ \boldsymbol{z}^{k}\!=\!col\{ z^{k}_{i} \}^{N}_{i=1}$ and analogously, 
	$$ \bar{x}^{k}= \frac{1}{N}\sum^{N}_{i=1}x_{i}^{k},\quad \overline{\boldsymbol{x}}^{k}=col\{ \overline{x}^{k} \}^{N}_{i=1},$$ 
	$$ \boldsymbol{s}^{k}\!=col\{ v^{k}_{i}-x^{k}_{i} \}^{N}_{i=1},\quad \overline{\boldsymbol{s}}^{k}=col\{ \frac{1}{N}\! \sum_{i=1}^{N} v_{i}^{k}-\bar{x}^{k} \}^{N}_{i=1}.$$ 
	Also, denote the averages global gradients by 
	$$ \overline{\nabla^{k} \!F}\!=\frac{1}{N}\! \sum_{i=1}^{N}\! \nabla\! f_{i}(x_{i}^{k}),\quad \overline{\nabla^{k}\! \boldsymbol{F}}\!=\!col\{\overline{\nabla^{k}\! F}\}_{i=1}^{N}.$$  
	Set $ \lambda(\mathcal{A})$ as the second largest eigenvalue of the adjcency matrix $ \mathcal{A} $. Select $ k_{0} $ as the smallest positive integer such that 
	\begin{equation}\label{ee3}
	\lambda(\mathcal{A}) \leq \frac{(k_{0}+1)^{2}-4(k_{0}+2)}{(k_{0}+2)(k_{0}+1)}.
	\end{equation}
	The existence of $ k_{0} $ is guaranteed because the adjacency matrix $\mathcal A$ is  symmetric and doubly stochastic \cite{wai2017decentralized}. Then we have the following result on the convergence rate of discrete algorithm \eqref{discrete-1}.

\begin{theorem} 
	Under \textit{Assumption 1}, set $ \delta=1 $, $ \beta^{k}= \frac{2}{k+1}$, and $ \bar{d} $ as the diameter  of $ \Omega $. For any $ k \geq 1$, there exist constants 
	$$ C_{x}=  \frac{1}{2}\sqrt{N}\bar{d}(k_{0}+1),$$ and $$C_{z} =\kappa \sqrt{N} \left(2 C_{x}+\bar{d}\right)   (k_{0}+1),$$ such that
	\begin{enumerate}[i).]
		\item all decision variables $x_i^k$ achieve consensus, \textit{i.e.},
		$$	\|\boldsymbol{x}^{k}-\overline{\boldsymbol{x}}^{k}\| \leq  \frac{2C_{x}}{k+1};
		$$
		\item 	all variable $z_i^k$ track the global gradient, \textit{i.e.}, 
		$$
		\|\boldsymbol{z}^{k}- \overline{\nabla^{k} \boldsymbol{F}} \|\leq \frac{2C_{z}}{k+1};
		$$
		\item 	the average $ \bar{x}^{k} $ converges to an optimal solution $ x^{*} $ with a  rate of $ O(1/k)$,  \textit{i.e.},
		$$
		F(\bar{x}^{k})-F\left(x^{*}\right)\leq \frac{8\bar{d}C_{z}+8\kappa\bar{d}C_{x} + 2\kappa \bar{d}^{2}  }{k+1}.$$
	\end{enumerate}
\end{theorem}	
\noindent\textbf{Proof.}

 i). Recall the compactness of $\Omega$ and this is naturally true for $k\in[1,k_{0}]$. For induction, assume 
 $$ \|\boldsymbol{x}^{k}-\overline{\boldsymbol{x}}^{k}\| \leq  \frac{2C_{x}}{k+1}, $$ for some $ k\geq k_{0} $. Then, we can verify that
$$
\begin{aligned}
\|\boldsymbol{x}^{k+1}-\bar{\boldsymbol{x}}^{k+1}\| 
&= \| \boldsymbol{A} \boldsymbol{x}^{k}+ \frac{2\boldsymbol{s}^{k}}{k+1}-(\overline{\boldsymbol{x}}^{k}+\frac{2\overline{\boldsymbol{s}}^{k}}{k+1}) \|\\
&\leq|\lambda(\mathcal{A})|\| \boldsymbol{x}^{k}-\bar{\boldsymbol{x}}^{k}\|+\frac{4\sqrt{N} \bar{d}}{k+1}  \\
&\leq \frac{2}{k+1} (|\lambda(\mathcal{A})|C_{x}+2\sqrt{N} \bar{d}) \\
&\leq    \frac{2C_{x}}{k+2},
\end{aligned}
$$
where  $\boldsymbol{A}=\mathcal{A} \otimes I_{n} $ and the last inequality follows from (\ref{ee3}). Thus, we have the conclusion according to induction. 

\noindent ii). Similarly, we can verify this for $k\in[1,k_0]$. For induction, assume 
$$
\|\boldsymbol{z}^{k}- \overline{\nabla^{k} \boldsymbol{F}} \|\leq \frac{2C_{z}}{k+1},
$$
for some $k\geq 0$. Then we can verify that
\begin{align}
&\|\boldsymbol{z}^{k+1}- \overline{\nabla^{k+1} \boldsymbol{F}}\|\nonumber\\
=&\|\boldsymbol{A}\boldsymbol{z}^{k}+G(\bm x^{k+1})-G(\bm x^{k})- \overline{\nabla^{k+1} \boldsymbol{F}}+ \overline{\nabla^{k} \boldsymbol{F}}- \overline{\nabla^{k} \boldsymbol{F}}\|\nonumber\\
&\leq\|\boldsymbol{A}\boldsymbol{z}^{k}- \overline{\nabla^{k} \boldsymbol{F}}\|+\|+G(\bm x^{k+1})-G(\bm x^{k})\|\nonumber\\
&+\| \overline{\nabla^{k+1} \boldsymbol{F}}- \overline{\nabla^{k} \boldsymbol{F}}\|.\label{123456}
\end{align}
The firsr term in \eqref{123456} satisfies 
\begin{align*}
\|\boldsymbol{A}\boldsymbol{z}^{k}- \overline{\nabla^{k} \boldsymbol{F}}\|\leq&\|\lambda(\mathcal{A})\|\|\boldsymbol{z}^{k}- \overline{\nabla^{k} \boldsymbol{F}}\|\\
\leq&\|\lambda(\mathcal{A})\|\frac{2C_{z}}{k+1}.
\end{align*}
The second term in \eqref{123456} satisfies
\begin{align*}
\|G(\bm x^{k+1})-G(\bm x^{k})\|=&\sqrt{\sum_{i=1}^N\|\nabla f_i(x_i^{k+1})-\nabla f_i(x_i^{k})\|^2}\\
\leq& \kappa\sqrt{\sum_{i=1}^N\|x_i^{k+1}-x_i^k\|^2}.
\end{align*}
Note that
\begin{align*}
&\|x_i^{k+1}-x_i^k\|\\
\leq&\|\sum_{j=1}^N a_{ij}(x_j^k-x_i^k)+\frac{2}{k+1}(v_i^k-x_i^k)\|\\
\leq&\sum_{j=1}^Na_{ij}(\|x_j^k-\bar x^k\|+\|x_i^k-\bar x^k\|+\frac{2}{k+1}\|v_i^k-x_i^k\|)\\
\leq&\frac{2(2C_x+\bar d)}{k+1}.
\end{align*}
So far, we have 
\begin{align*}
\|G(\bm x^{k+1})-G(\bm x^{k})\|\leq 2\kappa\sqrt{N}\frac{2C_x+\bar d}{k+1}.
\end{align*}
Next, the third term in \eqref{123456} satisfies 
\begin{align*}
\| \overline{\nabla^{k+1} \boldsymbol{F}}- \overline{\nabla^{k} \boldsymbol{F}}\|\leq&\|G(\bm x^{k+1}-G(\bm x^{k})\|\\
\leq&2\kappa\sqrt{N}\frac{2C_x+\bar d}{k+1}.
\end{align*}
To sum up, we obtain that
\begin{align*}
&\|\boldsymbol{z}^{k+1}- \overline{\nabla^{k+1} \boldsymbol{F}}\|\\
\leq& \|\lambda(\mathcal{A})\|\frac{2C_{z}}{k+1}+4\kappa\sqrt{N}\frac{2C_x+\bar d}{k+1}\\
\leq & \frac{2C_z}{k+2},
\end{align*}
where the last inequality follows from \eqref{ee3} and the definitions of constants $C_x$ and $C_z$. Hence, we get the conclusion according to induction.

\noindent iii).
Take 
$$g^{k}\triangleq F(\bar{x}^{k})-F\left(x^{*}\right).$$ It follows from \eqref{discrete-1} that 
$$
\bar{x}^{k+1}=\bar{x}^{k}+\frac{\beta^{k}}{N} \sum_{i=1}^{N} (v_{i}^{k}-\bar{x}^{k}).
$$
Since $F$ is $\kappa  $-Lipschitz  and $\Omega$ is compact, 
\begin{align}\label{f13}
g^{k+1} \leq g^{k}+\frac{\beta^{k}}{N} \sum_{i=1}^{N}\langle v_{i}^{k}-\bar{x}^{k}, \nabla F(\bar{x}^{k})\rangle+\frac{\kappa \bar{d}^{2}}{2}(\beta^{k})^{2}.
\end{align}
For  $i\in\mathcal{V}$,  because $ v_{i}^{k} \in \operatorname{argmin}_{v\in \Omega}\langle v,  z_{i}^{k}\rangle$, 
$$
\begin{aligned}
&\langle v_{i}^{k}-\bar{x}^{k}, \nabla F(\bar{x}^{k})\rangle\\ \leq&\langle v-\bar{x}^{k}, z_{i}^{k}\rangle+\bar{d}\|z_{i}^{k}-\nabla F(\bar{x}^{k})\|\\
\leq&\langle v-\bar{x}^{k}, \nabla F(\bar{x}^{k})\rangle+2 \bar{d}\|z_{i}^{k}-\nabla F(\bar{x}^{k})\|,\quad \forall v \in \Omega.
\end{aligned}
$$
Thus, take 
$$ v^{\prime}\in\operatorname{argmin}_{v\in \Omega} \langle v,  \nabla F(\bar{x}^{k}) \rangle .$$ 
Together with the optimality of $ v^{\prime} $ and the convexity of $ F $, 
\begin{align}\label{f14}
\langle v^{\prime}-\bar{x}^{k}, \nabla F(\bar{x}^{k}) \rangle &\leq\langle  x^{*}-\bar{x}^{k}, \nabla F(\bar{x}^{k})\rangle\nonumber\\ &\leq-g^{k}.
\end{align}
In addition,
\begin{align} \label{f12}
&\|z_{i}^{k}-\nabla F(\bar{x}^{k})\|\nonumber\\ \leq&\|z_{i}^{k}-\overline{\nabla^{k} F}\|+\|\overline{\nabla^{k} F}-\nabla F(\bar{x}^{k})\|\nonumber \\
\leq& \|z_{i}^{k}-\overline{\nabla^{k} F}\|+\frac{\kappa }{N} \sum_{i=1}^{N}\|x_{i}^{k}-\bar{x}^{k}\|.
\end{align} 	
When $ \beta^{k}= 2/(k+ 1 ) $,
it follows from the obtained conclusions in i) and ii) that, \eqref{f12} becomes
\begin{align}\label{f12.1}
&\|z_{i}^{k}-\nabla F(\bar{x}^{k})\| 
\leq  \frac{ 2C_{z} }{k+1}+\kappa  \frac{ 2C_{x} }{k+1}.
\end{align} 
By substituting (\ref{f14}) and \eqref{f12.1} into (\ref{f13}), we further have
\begin{align*}
&g^{k+1} \leq g^{k}-\frac{2g^{k}}{k+1}+\frac{8 \bar{d}\left(C_{z}+\kappa C_{x}\right)}{(k+1)^2}+\frac{2\kappa \bar{d}^{2}}{(k+1)^2},
\end{align*}
which yields 
$$ g^{k+1}\leq \frac{k-1}{k+1}g^{k} + O(1/k^{2}) .$$
Consequently,  by \cite[Lemma 4]{rivet1987introduction}, we have $ O(1/k) $ convergence rate for $ g^{k} $, that is, 
\begin{equation*}
g^{k} \leq  (8\bar{d}\left(C_{z}+\kappa C_{x}\right)+2\kappa \bar{d}^{2})/(k+1), \quad \forall\, k \geq 1.
\end{equation*}
This completes the proof.
\hfill$\square$

\section{Numerical Examples}\label{example}
\begin{figure}
	\begin{center}
		\includegraphics[height=6.5cm]{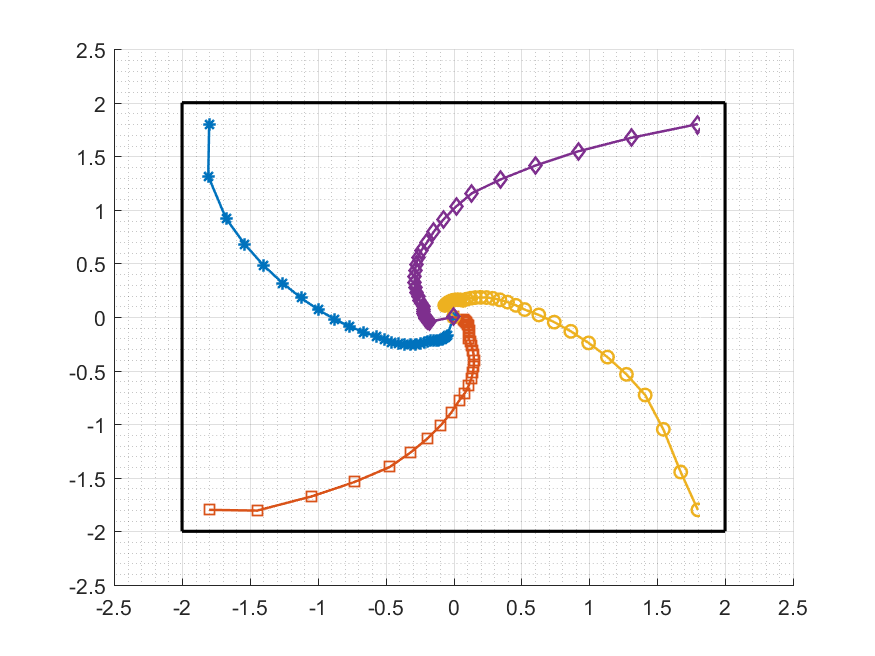}    
		\caption{trajectories of  the four agents' decision variables.}  
		\label{fig1}                                 
	\end{center}                                 
\end{figure}

\begin{table*}
	\caption{the average real running time of solving subproblems.}
	\renewcommand\arraystretch{2} 
	\begin{center}
		\setlength{\tabcolsep}{8mm}{
			\begin{tabular}{c|c|c|c|c|c}
				\hline
				\hline
				
				dimensions& n=16&n=64&n=256&n=1024&n=4096 \\
				\hline
				DCPf (msec)& 11.5 & 12.0  & 12.6 & 13.3 &14.1\\
				\hline
				DeFW (msec)& 11.5 & 12.1  & 12.8 & 13.1 &14.4\\
				\hline
				DCPb-Zeng (msec)&  8.7  & 13.2 &19.8 & 27.9&40.7\\
				\hline
				DCPb-Yang (msec)&  8.7  & 13.1 &19.5 & 28.4&43.0\\
				\hline
				\hline
		\end{tabular}}
	\end{center}
	\label{runtime}
\end{table*}

\begin{figure}

	\centering
	\subfloat[$n=16$]{
		\begin{minipage}{8cm} 
			\includegraphics[width=\textwidth]{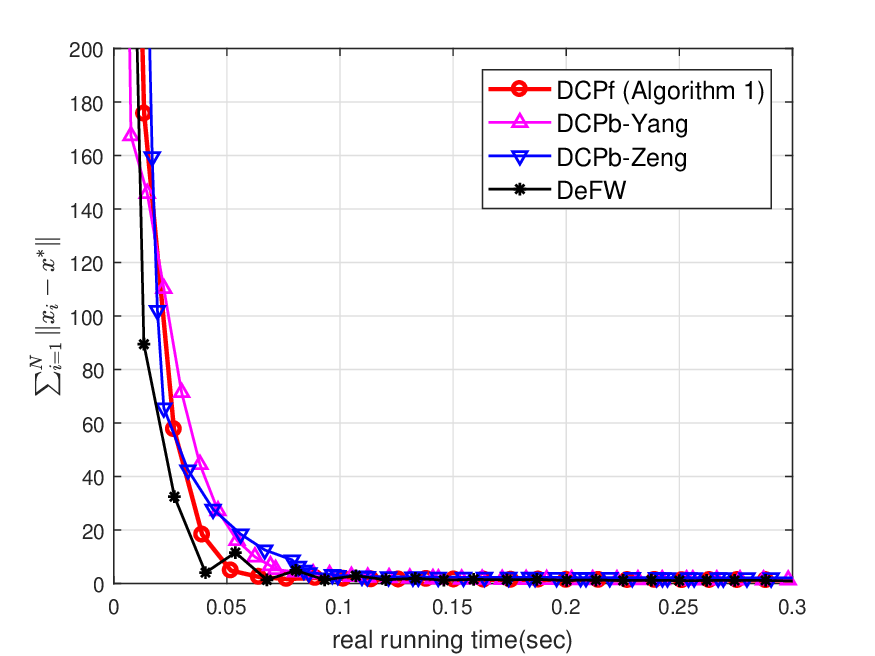} \\
		\end{minipage}
	}
	
	\subfloat[$n=64$]{
		\begin{minipage}{8cm}
			\includegraphics[width=\textwidth]{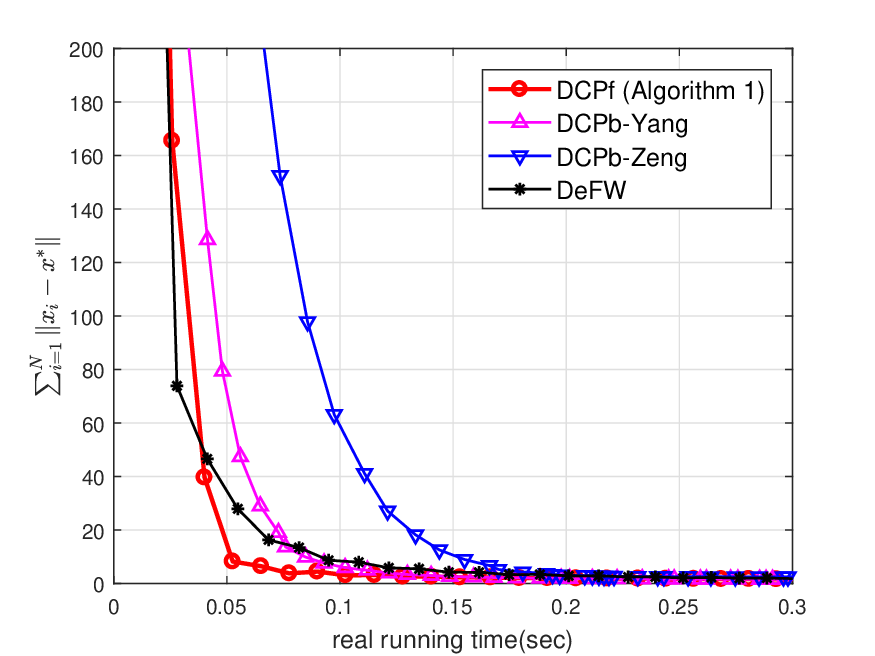} \\
			
		\end{minipage}
	}
	
	\subfloat[$n=256$]{
		\begin{minipage}{8cm}
			\includegraphics[width=\textwidth]{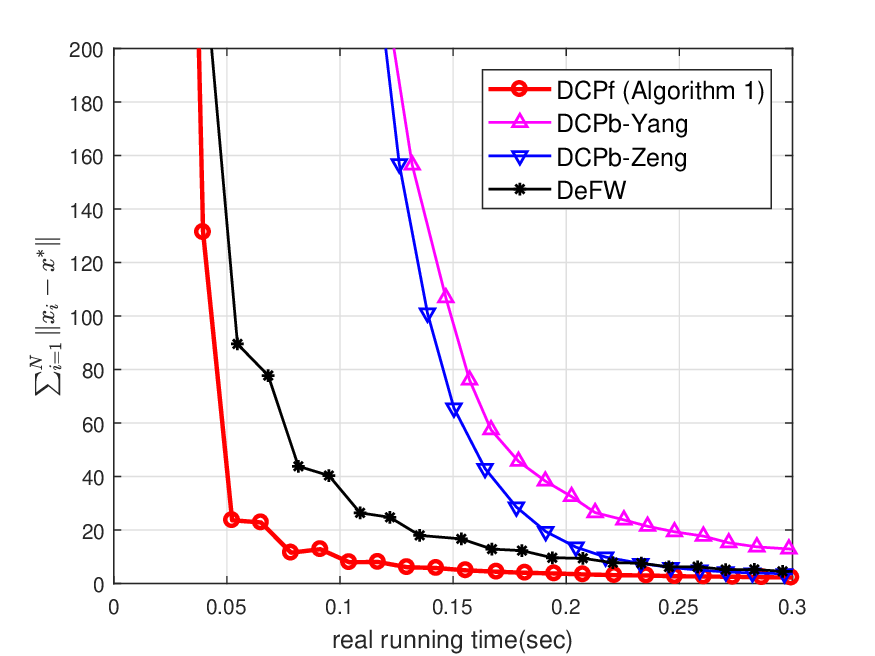} \\
			
		\end{minipage}
	}
	\caption{Optimal solution errors with different dimensions $n=16,64,256$.}
	\label{fig2}
\end{figure}

In this section, we carry on some experiments to illustrate the effectiveness of our approach in this paper.
We first take a simple example with only $N = 4$ agents and only $n=2$ dimensions of the decision variables to illustrate the trajectories of \textit{Algorithm \ref{alg:1}}. Set the local cost functions as
$$f_j(x)=(x_{j}^1-\frac{5}{3}+\frac{2}{3}j)^2+(x_{j}^2-\frac{5}{3}+\frac{2}{3}j)^2,\quad j=1,2,3,4.$$
The feasible sets is  
$$\Omega=\{x\in\mathbb R^2:-2\leq x^1\leq 2,-2\leq x^2\leq 2\}.$$
The initial locations are as follows,
\begin{align*}
x_1(0)&=\begin{bmatrix} -1.8 \\ 1.8 \end{bmatrix},\, x_2(0)=\begin{bmatrix} -1.8 \\ -1.8 \end{bmatrix},\\
x_3(0)&=\begin{bmatrix} 1.8 \\ 1.8 \end{bmatrix},\quad x_4(0)=\begin{bmatrix} 1.8 \\ -1.8 \end{bmatrix}.
\end{align*}
Under this circumstance, the global optimal solution should be the origin on Euclidean space.  We employ a directed ring graph as the communication network. Fig.\ref{fig1} shows the trajectory in the plane. The boundaries of the feasible set are in black, while the trajectories of  the four agents' decision variables are showed with different colors.

Next, we show the effectiveness of our distributed projection-free dynamics by comparisons. The number of agents is increased to $N=20$. As shown in \cite{jaggi2013revisiting}, the FW method works better than projection-based algorithms when  the vertices on the boundary of the constraint set are easily to find in high-dimensional decision spaces. Thus, along with quadratic local cost functions, we set the constraint set as $\|x\|_\infty\leq 2$. Then we choose different dimensions of decision variables and compare our distributed continuous-time projection-free algorithm (DCPf) with two distributed continuous-time projection-based algorithms, (DCPb-Yang) by \cite{yang2016multi}, (DCPb-Zeng) by \cite{zeng2016distributed}, and a decentralized discrete-time FW  algorithm (DeFW) by \cite{wai2017decentralized}. Since all these distributed algorithms are suitable for undirected graphs, we set an undirected ring graph as the communication network.

We set the dimensions of decision variables as the power of two. In Fig.\ref{fig2}, the $x$-axis is for the real running time (CPU time) in seconds, while the $y$-axis is for the optimal solution errors in each algorithm. We learn from Fig.\ref{fig2} that as the dimension increases, the real running time (CPU time) of projection-based algorithms is obviously longer than projection-free ones, because searching the vertices on the boundary of high-dimensional constraint sets (to solve a linear program) is faster than calculating a projection on high-dimensional constraint sets (to solve a quadratic program). Moreover, we can observe from Fig.\ref{fig2} that our DCPf is not second to DeFW over connected undirected graphs.

Furthermore, in Tab.\ref{runtime}, we list the average real running time of solving subproblems, \textit{i.e.}, linear programs or quadratic programs.
When the dimension is low, solving linear programs may take more time than solving  quadratic programs over such constraint sets. However, as the dimension increases explosively, solving quadratic programs in such situation turns to be difficult, but the time of solving linear programs still remains almost the same. That conforms with the advantage of projection-free approaches.

\section{Discussions}\label{conclusion}
In this paper, we developed a novel projection-free dynamics for solving distributed optimization with constraints. By employing the Frank-Wolfe method, we found a feasible descent direction by solving  optimization with a linear objective, which avoided solving high-dimensional subproblems in projection-based algorithms. In the distributed dynamics, we simultaneously made the decision variables achieve consensus and the local auxiliary variables track the global gradient. Then we gave an analysis on the convergence of the projection-free dynamics by the Lyapunov theory. As for the induced discrete-time format, we proved the convergence rate. Finally, we provided comparative illustrations to show the efficiency of our algorithm.

In the future, we may focus on some potential challenges to improve the result in this paper. First, whether fixed constants or self-adaptive settings are suitable for parameter $\beta(t)$ in the algorithm. Second, more advanced approaches should be investigated for directly analyzing the convergence rate of the continuous-time algorithm itself. Third, whether more general formulations of distributed constrained optimization can adopt analogous projection-free algorithms, such as resource allocation problems with coupled constraints. At last, as the advantages from the FW method itself, such projecion-free approaches fit a class of constrained optimization well, including polyhedrons, unit simplices, and unit squares. However, this approach does not always mean efficient for all kinds of constraints, and how to adopt similar ideas to handle diverse constraints is still a challenging but interesting task. Maybe mix-projection schemes for distributed optimization with different constraints will be developed in the near future.

\section*{Appendix}
{\appendices
	\section*{Proof of \textit{Lemma \ref{lemma1}}}
	For a convex set $C\subseteq\mathbb R^n$ and $x\in C$,  denote the normal cone to $C$ at $x$ by
$$\mathcal N_C(x)=\big\{v\in\mathbb R^{n}: v^{\rm T}(y-x)\leq 0,\quad\forall y\in C\big\},$$
and the tangent cone to  $C$ at $x$ by
$$\mathcal T_C(x)=\{\lim_{k\rightarrow\infty}\frac{x_k-x}{t_k}:x_k\in C,t_k>0,x_k\rightarrow x,t_k\rightarrow 0\}.$$
Let $P_{{\Omega}}(x_i(t))$ as the projection on $\Omega$ at point $x_i(t)$, which yields that $$x_i(t)-P_{\Omega}(x_i(t))\in\mathcal N_{\Omega}(x_i(t)).$$
Consider  $x_i(t)\in\Omega$ for $i\in\mathcal V$ and some $t\geq 0$.
Since $v_i(t)$, for $i=1,\dots,N$, are also selected from $\Omega$,  we have  
$$v_i(t)-x_i(t)\in\mathcal T_{\Omega}(x_i(t)),$$ 
and similarly, $$x_j(t)-x_i(t)\in\mathcal T_{\Omega}(x_i(t)).$$ Hence, it follows from the dynamics of \textit{Algorithm \ref{alg:1}} that
\begin{align*}\dot x_i(t)=\sum_{j=1}^{N}a_{ij}(x_j(t)\!-\!x_i(t))+\beta(t)(v_i(t)\!-\!x_i(t)),
\end{align*} which leads to $\dot x_i(t)
\in\mathcal T_{\Omega}(x_i(t))$.
On this basis, consider the energy function as $$E(t)=\frac{1}{2}\|x_i(t)-P_{{\Omega}}(x_i(t))\|^2.$$
Its derivative along the dynamics of \textit{Algorithm \ref{alg:1}} is
\begin{align*}\dot E(t)=\langle x_i(t)-P_{{\Omega}}(x_i(t)),\,\dot x_i(t)\rangle\leq0,
\end{align*}
where the last inequality holds because normal cones and tangent cones are orthogonal. Since $x_i(0)\in\Omega$, we have $E(0)=0$. Together with $E(t)\geq 0$ and $\dot E(t)\leq 0$, we have $E(t)\equiv 0$ at any time $t\geq 0$. This reveals that once all variables are located within ${{\Omega}}$ for some $t\geq 0$, they will not escape. Therefore, recalling  $x_i(0)\in\Omega$ for $i\in\mathcal V$, we complete the proof. \hfill $\square$

	\section*{Proof of \textit{Lemma \ref{lemma2}}}
Let us denote the funciton
$$h(t)=\exp\int_{0}^{t}\gamma(\tau)d\tau,$$ which implies  $$\lim_{t\rightarrow \infty}h(t)=\infty,\quad \dot h(t)=\gamma(t)h(t).$$ Recall that $s(t)\geq 0$ and 
$$\dot s(t)\leq -\gamma(t)s(t)+\gamma(t)\varepsilon(t).$$ After multiplying the both sides of the inequality above by $h(t)$, we can further obtain that
$$\frac{d}{dt}\big (s(t)h(t)\big)\leq \gamma(t)h(t)\varepsilon(t).$$
Then we integrate the above on the segment $(0,t)$ by the \textit{Comparison Lemma} in \cite{khalil2002nonlinear}, which leads to the following inequality
$$s(t)\leq\frac{s(0)}{h(t)}+\frac{1}{h(t)}\int_{0}^{t}\gamma(\tau)h(\tau)\varepsilon(\tau)d\tau.$$
Hereupon, we can give a detailed discussion based on the above result. On the one hand, if $$\int_{0}^{\infty}\gamma(\tau)h(\tau)\varepsilon(\tau)d\tau<\infty,$$ then $\lim_{t\rightarrow \infty}s(t)=0$ eventually.
On the other hand, if 
$$\int_{0}^{\infty}\gamma(\tau)h(\tau)\varepsilon(\tau)d\tau=\infty,$$
then it follows from L' Hospital rule that
\begin{align*}
\lim_{t\rightarrow \infty}\sup s(t)\leq&\lim_{t\rightarrow \infty}\frac{\gamma(t)h(t)\varepsilon(t)}{\gamma(t)h(t)}\\
&\quad\quad=\lim_{t\rightarrow \infty}\varepsilon(t)=0,
\end{align*}
which completes the proof. \hfill $\square$

}
\bibliographystyle{IEEEtran}        
\bibliography{ref}

\begin{thebibliography}{10}
\providecommand{\url}[1]{#1}
\csname url@samestyle\endcsname
\providecommand{\newblock}{\relax}
\providecommand{\bibinfo}[2]{#2}
\providecommand{\BIBentrySTDinterwordspacing}{\spaceskip=0pt\relax}
\providecommand{\BIBentryALTinterwordstretchfactor}{4}
\providecommand{\BIBentryALTinterwordspacing}{\spaceskip=\fontdimen2\font plus
\BIBentryALTinterwordstretchfactor\fontdimen3\font minus
  \fontdimen4\font\relax}
\providecommand{\BIBforeignlanguage}[2]{{%
\expandafter\ifx\csname l@#1\endcsname\relax
\typeout{** WARNING: IEEEtran.bst: No hyphenation pattern has been}%
\typeout{** loaded for the language `#1'. Using the pattern for}%
\typeout{** the default language instead.}%
\else
\language=\csname l@#1\endcsname
\fi
#2}}
\providecommand{\BIBdecl}{\relax}
\BIBdecl

\bibitem{nedic2010constrained}
A.~Nedic, A.~Ozdaglar, and P.~A. Parrilo, ``Constrained consensus and
  optimization in multi-agent networks,'' \emph{IEEE Transactions on Automatic
  Control}, vol.~55, no.~4, pp. 922--938, 2010.

\bibitem{yuan2012distributed}
D.~Yuan, S.~Xu, H.~Zhao, and L.~Rong, ``Distributed dual averaging method for
  multi-agent optimization with quantized communication,'' \emph{Systems \&
  Control Letters}, vol.~61, no.~11, pp. 1053--1061, 2012.

\bibitem{li2020distributed}
X.~Li, G.~Feng, and L.~Xie, ``Distributed proximal algorithms for multi-agent
  optimization with coupled inequality constraints,'' \emph{IEEE Transactions
  on Automatic Control}, 2020.

\bibitem{yuan2015regularized}
D.~Yuan, D.~W. Ho, and S.~Xu, ``Regularized primal--dual subgradient method for
  distributed constrained optimization,'' \emph{IEEE transactions on
  Cybernetics}, vol.~46, no.~9, pp. 2109--2118, 2015.

\bibitem{wang2015distributed}
X.~Wang, Y.~Hong, and H.~Ji, ``Distributed optimization for a class of
  nonlinear multiagent systems with disturbance rejection,'' \emph{IEEE
  transactions on Cybernetics}, vol.~46, no.~7, pp. 1655--1666, 2015.

\bibitem{yang2016multi}
S.~Yang, Q.~Liu, and J.~Wang, ``A multi-agent system with a
  proportional-integral protocol for distributed constrained optimization,''
  \emph{IEEE Transactions on Automatic Control}, vol.~62, no.~7, pp.
  3461--3467, 2016.

\bibitem{zeng2016distributed}
X.~Zeng, P.~Yi, and Y.~Hong, ``Distributed continuous-time algorithm for
  constrained convex optimizations via nonsmooth analysis approach,''
  \emph{IEEE Transactions on Automatic Control}, vol.~62, no.~10, pp.
  5227--5233, 2016.

\bibitem{yi2018distributed}
X.~Yi, L.~Yao, T.~Yang, J.~George, and K.~H. Johansson, ``Distributed
  optimization for second-order multi-agent systems with dynamic
  event-triggered communication,'' in \emph{Conference on Decision and
  Control}.\hskip 1em plus 0.5em minus 0.4em\relax IEEE, 2018, pp. 3397--3402.

\bibitem{shi2018finite}
X.~Shi, J.~Cao, G.~Wen, and M.~Perc, ``Finite-time consensus of opinion
  dynamics and its applications to distributed optimization over digraph,''
  \emph{IEEE transactions on Cybernetics}, vol.~49, no.~10, pp. 3767--3779,
  2018.

\bibitem{lu2018distributed}
K.~Lu, G.~Jing, and L.~Wang, ``Distributed algorithms for searching generalized
  nash equilibrium of noncooperative games,'' \emph{IEEE transactions on
  Cybernetics}, vol.~49, no.~6, pp. 2362--2371, 2018.

\bibitem{khalil2002nonlinear}
H.~K. Khalil and J.~W. Grizzle, \emph{Nonlinear {S}ystems}.\hskip 1em plus
  0.5em minus 0.4em\relax Prentice {H}all Upper Saddle River, NJ, 2002, vol.~3.

\bibitem{ruszczynski2011nonlinear}
A.~Ruszczynski, \emph{Nonlinear {O}ptimization}.\hskip 1em plus 0.5em minus
  0.4em\relax Princeton {U}niversity {P}ress, 2011.

\bibitem{yi2016initialization}
P.~Yi, Y.~Hong, and F.~Liu, ``Initialization-free distributed algorithms for
  optimal resource allocation with feasibility constraints and application to
  economic dispatch of power systems,'' \emph{Automatica}, vol.~74, pp.
  259--269, 2016.

\bibitem{zhu2018projected}
Y.~Zhu, W.~Yu, G.~Wen, and G.~Chen, ``Projected primal--dual dynamics for
  distributed constrained nonsmooth convex optimization,'' \emph{IEEE
  Transactions on Cybernetics}, vol.~50, no.~4, pp. 1776--1782, 2018.

\bibitem{yang2019survey}
T.~Yang, X.~Yi, J.~Wu, Y.~Yuan, D.~Wu, Z.~Meng, Y.~Hong, H.~Wang, Z.~Lin, and
  K.~H. Johansson, ``A survey of distributed optimization,'' \emph{Annual
  Reviews in Control}, vol.~47, pp. 278--305, 2019.

\bibitem{liang2020distributed}
S.~Liang, X.~Zeng, G.~Chen, and Y.~Hong, ``Distributed sub-optimal resource
  allocation via a projected form of singular perturbation,''
  \emph{Automatica}, vol. 121, p. 109180, 2020.

\bibitem{frank1956algorithm}
M.~Frank, P.~Wolfe \emph{et~al.}, ``An algorithm for quadratic programming,''
  \emph{Naval {R}esearch {L}ogistics {Q}uarterly}, vol.~3, no. 1-2, pp.
  95--110, 1956.

\bibitem{chen123distributed}
G.~Chen, Y.~Ming, Y.~Hong, and P.~Yi, ``Distributed algorithm for
  $\varepsilon$-generalized nash equilibria with uncertain coupled
  constraints,'' \emph{Automatica}, vol. 123, p. 109313.

\bibitem{jaggi2013revisiting}
M.~Jaggi, ``Revisiting {F}rank-{W}olfe: Projection-free sparse convex
  optimization,'' in \emph{International {C}onference on {M}achine {L}earning},
  2013, pp. 427--435.

\bibitem{jiang2022distributed}
X.~Jiang, X.~Zeng, L.~Xie, J.~Sun, and J.~Chen, ``Distributed stochastic
  projection-free solver for constrained optimization,'' \emph{arXiv preprint
  arXiv:2204.10605}, 2022.

\bibitem{cui2014generalized}
Z.~Cui, H.~Chang, S.~Shan, and X.~Chen, ``Generalized unsupervised manifold
  alignment,'' \emph{Advances in Neural Information Processing Systems},
  vol.~27, pp. 2429--2437, 2014.

\bibitem{chapel2020partial}
L.~Chapel, M.~Z. Alaya, and G.~Gasso, ``Partial optimal tranport with
  applications on positive-unlabeled learning,'' \emph{Advances in Neural
  Information Processing Systems}, vol.~33, pp. 2903--2913, 2020.

\bibitem{zhang2017bio}
X.~Zhang and S.~Mahadevan, ``A bio-inspired approach to traffic network
  equilibrium assignment problem,'' \emph{IEEE Transactions on Cybernetics},
  vol.~48, no.~4, pp. 1304--1315, 2017.

\bibitem{garber2015faster}
D.~Garber and E.~Hazan, ``Faster rates for the {F}rank-{W}olfe method over
  strongly-convex sets,'' in \emph{International {C}onference on {M}achine
  {L}earning}, 2015, pp. 541--549.

\bibitem{zhang2017projection}
W.~Zhang, P.~Zhao, W.~Zhu, S.~C. Hoi, and T.~Zhang, ``Projection-free
  distributed online learning in networks,'' in \emph{International
  {C}onference on {M}achine {L}earning}, 2017, pp. 4054--4062.

\bibitem{zhang2020quantized}
M.~Zhang, L.~Chen, A.~Mokhtari, H.~Hassani, and A.~Karbasi, ``Quantized
  {F}rank-{W}olfe: Faster optimization, lower communication, and projection
  free,'' in \emph{International Conference on Artificial Intelligence and
  Statistics}, 2020, pp. 3696--3706.

\bibitem{jacimovic1999continuous}
M.~Jacimovic and A.~Geary, ``A continuous conditional gradient method,''
  \emph{Yugoslav {J}ournal of {O}perations {R}esearch}, vol.~9, no.~2, pp.
  169--182, 1999.

\bibitem{wai2017decentralized}
H.-T. Wai, J.~Lafond, A.~Scaglione, and E.~Moulines, ``Decentralized
  {F}rank-{W}olfe algorithm for convex and nonconvex problems,'' \emph{IEEE
  Transactions on Automatic Control}, vol.~62, no.~11, pp. 5522--5537, 2017.

\bibitem{gharesifard2013distributed}
B.~Gharesifard and J.~Cort{\'e}s, ``Distributed continuous-time convex
  optimization on weight-balanced digraphs,'' \emph{IEEE Transactions on
  Automatic Control}, vol.~59, no.~3, pp. 781--786, 2013.

\bibitem{nedic2017achieving}
A.~Nedic, A.~Olshevsky, and W.~Shi, ``Achieving geometric convergence for
  distributed optimization over time-varying graphs,'' \emph{SIAM Journal on
  Optimization}, vol.~27, no.~4, pp. 2597--2633, 2017.

\bibitem{pu2020distributed}
S.~Pu and A.~Nedi{\'c}, ``Distributed stochastic gradient tracking methods,''
  \emph{Mathematical Programming}, pp. 1--49, 2020.

\bibitem{qiu2016distributed}
Z.~Qiu, S.~Liu, and L.~Xie, ``Distributed constrained optimal consensus of
  multi-agent systems,'' \emph{Automatica}, vol.~68, pp. 209--215, 2016.

\bibitem{horn2012matrix}
R.~A. Horn and C.~R. Johnson, \emph{Matrix {A}nalysis}.\hskip 1em plus 0.5em
  minus 0.4em\relax Cambridge {U}niversity {P}ress, 2012.

\bibitem{rivet1987introduction}
A.~Rivet and A.~Souloumiac, ``Introduction to optimization,'' in
  \emph{Optimization Software, Publications Division}.\hskip 1em plus 0.5em
  minus 0.4em\relax Citeseer, 1987.

\end{thebibliography}

\end{document}